\documentclass[12pt]{amsart}

\usepackage[letterpaper,top=3cm,bottom=3cm,left=2cm,right=2cm,marginparwidth=1.75cm]{geometry}
\usepackage{graphicx,xcolor}
\usepackage[colorlinks=true, allcolors=violet]{hyperref}

\usepackage{amssymb,amsfonts}
\usepackage[all,arc]{xy}
\usepackage{enumerate}
\usepackage{mathrsfs}
\usepackage{mathtools}
\usepackage{MnSymbol}

\newtheorem{thm}{Theorem}[section]
\newtheorem{cor}[thm]{Corollary}
\newtheorem{prop}[thm]{Proposition}
\newtheorem{lem}[thm]{Lemma}

\theoremstyle{definition}

\newtheorem{exmp}[thm]{Example}

\theoremstyle{remark}
\newtheorem{rem}[thm]{Remark}

\newcommand{\Z}{\mathbb{Z}}

\newcommand{\R}{\mathbb{R}}
\newcommand{\C}{\mathbb{C}}
\newcommand{\F}{\mathbb{F}}

\makeatletter
\let\c@equation\c@thm
\makeatother
\numberwithin{equation}{section}

\title{Continuous Inverse Ambiguous Functions on Lie Groups}

\author{David Schmitz}
\author{Sadman Rahman}
\author{Anthony Kindness}
\address{Department of Mathematics, North Central College, Naperville, IL.}

\begin{document}

\begin{abstract}
In Schmitz (Aequ Math 91:373 - 389, 2017), the first author defines an inverse ambiguous function on a group $G$ to be a bijective function $f : G \rightarrow G$ satisfying the functional equation $f^{-1}(x) = f(x^{-1})$ for all $x \in G$.  In this paper, we investigate the existence of continuous inverse ambiguous functions on various Lie groups.  In particular, we look at tori, elliptic curves over various fields, vector spaces, additive matrix groups, and multiplicative matrix groups.
\end{abstract}

\maketitle

\tableofcontents

\section{Introduction}
We define an inverse ambiguous function to be a bijective function $f : G \rightarrow G$, where $G$ is a group, that satisfies the functional equation 
\begin{equation}
    f^{-1}(x) = (f(x))^{-1}
\end{equation} 
for all $x \in G$. This condition is equivalent to
\begin{equation}
    f(f(x)) = x^{-1}
\end{equation}
for all $x \in G$, where $x^{-1}$ denotes the inverse of $x$ in $G$, which in turn implies that $f(x^{-1}) = f(x)^{-1}$ for all $x \in G$ ([5], Corollary 2.3).

In [5], the first author investigated the existence of inverse ambiguous functions on the additive and multiplicative groups associated to fields, with the added assumption of continuity for topological fields.  In [6], the existence of inverse ambiguous functions on some non-abelian groups was studied.  For finite groups, we have the following two facts [5, Theorem 3.3 and Theorem 4.1]:

\begin{thm}  \label{generalexistence} Let $G$ be a finite group, and let $S(G)$ denote the subset of $G$ consisting of all elements of order $1$ or $2$.  Then
\begin{enumerate}
    \item  If $f:G \rightarrow G$ is an inverse ambiguous function, then $f(S(G)) = S(G)$.
    \item  An inverse ambiguous function exists on $G$ if and only if the number of elements in $G - S(G)$ is a multiple of 4.
\end{enumerate}
\end{thm}

In this paper, we focus on the existence of continuous inverse ambiguous functions on topological groups, primarily Lie groups.  Recall that a topological group is a topological space endowed with a continuous binary operation satisfying the usual group operations, where the inverse map $x \mapsto x^{-1}$ is continuous.  We note that if $G$ is a topological group and $f: G \rightarrow G$ is a continuous inverse ambiguous function, then according to Equation (1), the inverse function $f^{-1}$ must also be continuous since it is the composition of two continuous functions.  In other words, a continuous inverse ambiguous function is necessarily a homeomorphism.  Therefore, if $G$ is a topological group and Homeo($G$) is the associated homeomorphism group (under composition), then the existence of a continuous inverse ambiguous function $f$ on $G$ is more or less equivalent to existence of a ``square root of $-1$", whereby $-1$ we mean the inversion map $x \mapsto x^{-1}$.  To put it a different way, $f$ generates a cyclic subgroup of order 4 in Homeo($G$). 

Using algebraic and topological methods, we will consider the question of existence of inverse ambiguous functions on the spherical groups $S^1$ and $S^3$, tori, elliptic curves (over various fields), and matrix groups. \\

\section{Tori}
We begin with the following lemma.

\begin{lem}  \label{circleno} There are no continuous inverse ambiguous functions on the circle group, $S^1$.
\end{lem}

\begin{proof}
By way of contradiction, assume that $f: S^1 \rightarrow S^1$ is an inverse ambiguous function.  By definition, $f \circ f = \iota$, where $\iota(z) = z^{-1}$ for all $z \in S^1$.  Since $f$ and $\iota$ are homeomorphisms on $G$, they induce group automorphisms $f_*$ and $\iota_*$, respectively, on the fundamental group $\pi_1(S^1) = \Z$ such that $f_* \circ f_* = \iota_*$.  Because $\iota$ reverses the direction of loops in $S^1$, it follows that $\iota_*(n) = -n$ for all $n \in \Z$.  However, the only group automorphisms on $\Z$ are the identity map and the map $n \mapsto -n$, and either of these composed with itself produces the identity map, not $\iota_*$.  We have reached a contradiction, and so there are no inverse ambiguous functions on $S^1$.
\end{proof}

The fundamental group can similarly be used to show the non-existence of inverse ambiguous functions on tori $\mathbb{T}^n$ when $n$ is odd.  For this, we recall the well-known result [1, Proposition 1.12] that for path-connected spaces $X_i$, $i = 1, \dots, n$, we have 
$$ \pi_1\left(\prod_{i=1}^n X_i\right) \cong \prod_{i=1}^n \pi_1(X_i) . $$

\begin{thm}  If $n$ is a positive odd integer, then there exist no continuous inverse ambiguous functions on the torus $\mathbb{T}^n$.
\end{thm}

\begin{proof}  By way of contradiction, assume that $f : \mathbb{T}^n \rightarrow \mathbb{T}^n$ is an inverse ambiguous function, where $n$ is odd.  As in the proof of the previous lemma, $f$ and $\iota$ are induce group automorphisms $f_*$ and $\iota_*$, respectively, on the fundamental group $\pi_1(\mathbb{T}^n) \cong \Z^n$ such that $f_* \circ f_* = \iota_*$. Under the isomorphism $\pi_1(\mathbb{T}^n) \cong \Z^n$, $\iota_*$ becomes the map $$\iota_*: (m_1, m_2,\dots, m_n) \mapsto (-m_1, -m_2, \dots, -m_n).$$  Thus, as an isomorphism of abelian groups (or $\Z$-modules), the map $f_*$ can be associated to an element in $A_f \in GL_n(\Z)$ such that $(A_f)^2 = -I_n$, where $I_n$ is the $n \times n$ identity matrix.  Taking determinants of both sides, we get $\det(A_f)^2 = (-1)^n = -1$, which is impossible since $\det(A_f)^2$ is an integer.  We have reached a contradiction, and so there are no inverse ambiguous functions on $\mathbb{T}^n$ when $n$ is odd.
\end{proof}

On the other hand, when $n$ is even, there do exist inverse ambiguous functions on $\mathbb{T}^n$.  We begin with looking at $\mathbb{T}^2$.  Define the function $f: S^1 \times S^1 \rightarrow S^1 \times S^1$ by
\[ f(z, w) = \left(w, z^{-1}\right). \]
Then, it is easy to verify that $f$ is a homeomorphism from $\mathbb{T}^2$ to itself such that 
$$ f(f(z,w)) = f\left(w,z^{-1}\right) = \left(z^{-1},w^{-1}\right) = (z,w)^{-1} $$
for all $z, w \in S^1$.  Hence, $f$ is inverse ambiguous.  To show the existence of inverse ambiguous functions on other tori $\mathbb{T}^n$ with $n$ even, we can the following, whose proof is straightforward.

\begin{prop}
Let $G$ and $H$ be topological groups and suppose that $g: G \rightarrow G$ and $h : H \rightarrow H$ are continuous inverse ambiguous functions.  Then $f = (g, h)$ is a continuous inverse ambiguous function defined on $G \times H$.
\end{prop}

Applying this proposition inductively, we obtain the following result.

\begin{thm}  If $n$ is a positive even integer, then there exist continuous inverse ambiguous functions on the torus $\mathbb{T}^n$.
\end{thm}

A specific example of a continuous (smooth, in fact) inverse ambiguous function on $\mathbb{T}^{2n}$ is defined by
\[ g:(z_1, w_1; z_2, w_2; \dots ; z_n, w_n) \mapsto \left(w_1, z_1^{-1}; w_2, z_2^{-1}; \dots ; w_n, z_n^{-1}\right). \]

\begin{rem} \label{evenproducts}
More generally, if $G$ is \textit{any} topological group and $H = G^{2n}$ (a direct product of an \textit{even} number of copies of $G$), then the function $g$ described above is a continuous inverse ambiguous function on $H$.
\end{rem}

\section{Elliptic Curves}

An elliptic curve $E$ over a field $K$ is a smooth 1-dimensional projective variety of genus 1.  Up to isomorphism, an elliptic curve consists of the set of solutions $[x, y, 1] \in \mathbb{P}^2(K)$ of a Weierstrass equation 
\[ Y^2 + a_1XY + a_3Y = X^3 + a_2X^2 + a_4X + a_6,\] with coefficients in $K$, together with the ``point at infinity" $O = [0,1,0]$.  Over fields of characteristic not equal to 2 or 3, a change of variables can transform the Weierstrass equation into the form 
\[Y^2 = X^3 + aX + B,\]
with $a, b \in K$.  Smoothness is equivalent to the discriminant $\Delta_E = -16(4a^3 + 27b^2)$ being non-zero. (See [7], Chapter 3.) Furthermore, it is not difficult to show (using differential calculus, for example) that $\Delta_E$ is positive (resp., negative), when $X^3 + aX + b$ has 1 real root (resp., 3 real roots). 

An elliptic curve over $\C$ is isomorphic (as a Lie group) to a complex torus $\C/\Lambda$, where $\Lambda$ is a lattice generated over $\Z$ by two $\R$-linearly independent complex numbers $\omega_1, \omega_2$:
\[\Lambda = \{ n_1 \omega_1 + n_2 \omega_2 \: | \: n_1, n_2 \in \Z \}.\]
(See [7, Corollary 5.1.1].)

\begin{thm} Every elliptic curve over $\C$ admits a continuous inverse ambiguous function.
\end{thm}

\begin{proof}  Using the notation above, the $\R$-linear independence of $\omega_1$ and $\omega_2$ implies that $(\omega_1, \omega_2)$ is an $\R$-basis of $\C$.  Therefore, we can define an $\R$-linear transformation $f: \C \rightarrow \C$ by $f(c_1\omega_1 + c_2\omega_2) = -c_2\omega_1 + c_1\omega_2$.  This function is a continuous isomorphism of $\R$-vector spaces with the property that 
\[f(f(c_1\omega_1 + c_2\omega_2)) = f(-c_2\omega_1 + c_1\omega_2) = -c_1\omega_1 - c_2\omega_2.\] 
In other words, $f(f(z)) = -z$ for all $z \in \C$. Moreover, $f(\Lambda) = \Lambda$.  Hence, $f$ induces a well-defined homeomorphism $\tilde{f} : \C/\Lambda \rightarrow \C/\Lambda$ where $\tilde{f}(z + \Lambda) = -z + \Lambda$.  In other words, $\tilde{f}$ is a continuous inverse ambiguous function on $\C/\Lambda$.
\end{proof}

It is known that the group of real points of an elliptic curve $E$ over $\R$ is isomorphic (as Lie groups) to either $S^1$ (which happens when $\Delta_E < 0$) or $S^1 \times \Z_2$ (which happens when $\Delta_E > 0$) [2, Introduction, Proposition 7.2].  In the next theorem, we show that a continuous inverse ambiguous function exists on the elliptic curve in the latter case, but not the former.

\begin{thm}  Let $E(\R)$ be the Lie group of real points on an elliptic curve over $\R$.  Then $E(\R)$ admits a continuous inverse ambiguous function if and only if $\Delta_E > 0$.
\end{thm}

\begin{proof}  First, suppose that $\Delta_E < 0$, then $E(\R) \cong S^1$.  According to Lemma~\ref{circleno}, no continuous inverse ambiguous functions can exist on $E(\R)$.  

Next, assume $\Delta_E > 0$.  Then $E(\R) \cong S^1 \times \Z_2$.  The two connected components of this space are the subgroup $S^1 \times \{0\}$ and its complement $S^1 \times \{1\}$.  Therefore, we can define a continuous function on $S^1 \times \Z_2$ that maps points of the form $(z, 0)$ to $(z, 1)$ and points of the form $(w, 1)$ to $(w^{-1},0)$.  One can easily check that $f(f(z,0)) = (z^{-1},0) = (z, 0)^{-1}$ and $f(f(z,1)) = (z^{-1},1) = (z, 1)^{-1}$ for every $z \in S^1$.  This shows that there exists a continuous inverse ambiguous function on $S^1 \times \Z_2$, and so the same is true for $E(\R)$.
\end{proof}

In fact, if $G$ is \textit{any} topological group, then the function defined by 
\[ f((z, 0)) = (z,1), \hspace{.4in} \textrm{and} \hspace{.4in} f((w,1)) = \left(w^{-1}, 0\right) \]
is a continuous inverse ambiguous function on $G \times \Z_2$.  With a little more effort, we can generalize this construction to some other finite groups $H$ besides $\Z_2$.

\begin{prop}  
    Let $G$ be any topological group and $H$ any finite group satisfying the following two conditions:
    \begin{enumerate}
        \item The number of self-invertible elements in $H$ is even.
        \item The number of non-self-invertible elements in $H$ is a multiple of 4. 
    \end{enumerate}
    Then $G \times H$ admits a continuous inverse ambiguous function.
\end{prop}

This proposition applies, in particular, to $A_4$, $A_5$, $A_n$ ($n \ge 8$), $S_n$ ($n \ge 6$) (see [6]), and cyclic groups $\Z_n$ with $n \equiv 2 \mod 4$ (see [5]).

\begin{proof}  
    With $G$ and $H$ satisfying the assumptions in the proposition, enumerate the self-invertible elements of $H$ as $$S(H) = \{a_1, a_2, ..., a_n, b_1, b_2, ..., b_n\}.$$  As in [5, Corollary 3.5], the non-self-invertible elements can be partitioned into subsets of the form $\{c_i, d_i, c_i^{-1}, d_i^{-1}\}$ for $i = 1, ..., m$.  Taking advantage of the fact that the subspaces \{$G \times \{x\}\}_{x \in H}$ are unions of connected components of $G \times H$, we can define a continuous function $f : G \times H \rightarrow G \times H$ by the following rules:
    \[\begin{array}{rlcl}
    f((g, a_i)) &= (g, b_i) & & \textrm{for } i = 1, \dots, n \\ 
    f((g, b_i)) &= (g^{-1}, a_i) & & \textrm{for } i = 1, \dots, n \\
    f((g, c_j)) &= (g, d_j) & & \textrm{for } j = 1, \dots, m \\
    f((g, d_j)) &= (g^{-1}, c_j^{-1}) & & \textrm{for } j = 1, \dots, m \\
    f((g, c_j^{-1})) &= (g, d_j^{-1}) & & \textrm{for } j = 1, \dots, m \\
    f((g, d_j^{-1})) &= (g^{-1}, c_j) & & \textrm{for } j = 1, \dots, m. 
    \end{array} \]
    It is a straightforward exercise to verify that $f(f((g,h))) = (g^{-1}, h^{-1})$ for all elements $(g, h)\in G \times H$.
\end{proof}

Finally, we consider elliptic curves defined over  finite fields.  For simplicity, we will assume $E$ is defined by a Weierstrass equation of the form 
\begin{equation}
    Y^2 = X^3 + aX + b
\end{equation} 
where $g(X) = X^3 + aX + b \in \Z[X]$.  If $\F_q$ ($q$ a power of a prime) is a finite field and $\Delta_E = -16(4a^3 + 27b^2)$ is non-zero in $\F_q$, then the group $E(\F_q)$ is finite, and with respect to the discrete topology, any inverse ambiguous function on $E(\F_q)$ is continuous. Therefore, according to Theorem~\ref{generalexistence}, the existence of an inverse ambiguous function on $E(\F_q)$ amounts to showing that the number of non-self-invertible elements in $E(\F_q)$ is a multiple of 4.

The self-invertible elements in $E(\F_q)$ are the point at infinity, $O$, and the elements of order 2.  The latter are known to be of the form $(r, 0)$, where $r$ is a root of $g(X)$ ([8], Section 2.1).  
Under the assumption that $\Delta_E \ne 0$ in $\F_q$, it follows that $g(X)$ has either 0, 1, or 3 roots in $\F_q$.  On the other hand, if $s \in \F_q$ is \textit{not} a root of $g(X)$, then the equation $Y = g(s)$ will have either 0 or 2 solutions for $Y$ in $\F_q$.  This implies that the number of non-self-invertible elements in $E(\F_q)$ will always be even, and it will be a multiple of 4 if and only if the set 
\[ N_q = \{ s \in \F_q \mid g(s) \textrm{ is a non-zero square in } \F_q \} \]
has an even number of elements.

Suppose $g(X)$ has 3 roots in $E(\F_q)$.  This would mean that $S(E(\F_q))$ is a subgroup of order 4 in $E(\F_q)$, which in turn implies, according to LaGrange's Theorem, that the order of $E(\F_q)$, and hence $N_q$, is a multiple of 4.  We have proved the following:

\begin{prop}  \label{elliptic3roots} Suppose $g(X) = X^3 + aX + b$ has three distinct roots in $\F_q$.  Then the elliptic curve defined by the Weierstrass equation $Y^2 = g(X)$ admits a (continuous) inverse ambiguous function.
\end{prop}

We now consider a couple of examples: 

\begin{exmp}  Consider the group of $\F_q$-points on the elliptic curve defined by the Weierstrass equation $Y^2 = X^3 - X$.  Then $\Delta_E = 64$, which is non-zero in $\F_q$ for every odd $q$.  For every such $q$, $g(X) = X^3 - X$ has three distinct roots in $E(\F_q)$, namely $-1, 0, 1$.  Therefore, by Proposition~\ref{elliptic3roots}, $E(\F_q)$ admits a (continuous) inverse ambiguous function for every odd $q$. 
\end{exmp}

\begin{exmp}  Consider the Weierstrass equation $Y^2 = X^3 + X$.  Then $\Delta_E = -64$, which is non-zero in $\F_q$ for every odd $q$.  For every such $q = p^n$ for which $p \equiv 1 \mod 4$ or $p \equiv 3 \mod 4$ with $n$ even, $g(X) = X^3 + X$ has three distinct roots in $E(\F_q)$, namely $0$ and the 2 square roots of $-1$.  For these values of $q$, Proposition~\ref{elliptic3roots} guarantees that $E(\F_q)$ admits a (continuous) inverse ambiguous function. Therefore, by Proposition~\ref{elliptic3roots}, $E(\F_q)$ admits a (continuous) inverse ambiguous function. 

However, if $q = p^n$ with $p \equiv 3 \mod 4$ and $n$ odd, then $g(X)$ has exactly 1 root ($X = 0$) in $\F_q$ because $p^n = X^2 + 1$ has no solution in $\Z$ ([4], Theorem 13.6).    Among the non-zero elements of $\F_q$ there are $(q-1)/2$ squares and $(q-1)/2$ non-squares, and the subset of non-squares is collectively the additive inverses of the subset of squares (because $-1$ is not a square).  Because $g(X)$ is an odd polynomial, exactly one of $g(a)$ and $g(-a)$ will be a square in $\F_q$ for each non-zero $a \in \F_q$.  Since $(q-1)/2$ is odd, it follows that the $|N_q| \equiv 2 \mod 4$, which means that there are no (continuous) inverse ambiguous functions defined on $E(\F_q)$. 
\end{exmp}

\section{Vector Spaces and Spheres}

We next look into the existence of continuous inverse ambiguous functions on vectors spaces.  Our investigation will consider vector spaces over various fields of scalars.  First, if $V$ is an $n$-dimensional complex vector space, then we see immediately that the linear transformation $f(v) = iv$ is a continuous inverse ambiguous function defined on $V$.

Next, let $V$ be an $n$-dimensional vector space over the finite field $\F_q$ where $q = p^r$ for some prime $p$.  Using the fact that the group $\F_q^\times$ is cyclic, it follows that if $q \equiv 1 \mod 4$ or if $p = 2$, then there exists $\alpha \in \F_q$ such that $\alpha^2 = -1$, and so $f : V \rightarrow V$ defined by $f(v) = \alpha v$ is a (continuous) inverse ambiguous function on $V$.  Now suppose $q \equiv 3 \mod 4$.  Then $V$ has $q^n$ elements, only one of which, $0_V$, is self-invertible.  Therefore, according to Theorem~\ref{generalexistence}, there exists a (continuous) inverse ambiguous function on $V$ if and only if $q^n \equiv 1 \mod 4$, which occurs if and only if $n$ is even. In summary:

\begin{thm}  Let $V$ be an $n$-dimensional vector space over the finite field $\F_q$, where $q$ is a power of a prime $p$. 
\begin{enumerate}
    \item  If $q$ is even or $q \equiv 1 \mod 4$, then $V$ admits a (continuous) inverse ambiguous function.
    \item  If $q \equiv 3 \mod 4$, then $V$ admits a (continuous) inverse ambiguous function if and only if $n$ is even.
\end{enumerate}
\end{thm}

We next turn our attention to finite-dimensional real vector spaces.  If $V = \R^{2n}$ is even-dimensional, then using Remark~\ref{evenproducts} it follows that

\begin{equation}  \label{imap}
f : (x_1, y_1; x_2, y_2; ...; x_n, y_n) \mapsto (y_1, -x_1; y_2, -x_2; ...; y_n, -x_n) 
\end{equation}

\vspace{5pt}

\noindent is a continuous inverse ambiguous function on $V$.  If $V$ is canonically identified with $\C^n$, then $f$ is the function that multiplies each coordinate by $i$.  

On the other hand, if $V = \R^{2n+1}$ is odd-dimensional, then we will demonstrate that there exists no continuous inverse ambiguous function on $V$.  On the contrary, suppose $f : V \rightarrow V$ were such a function.  Since $0$ is the only self-invertible element in $V$, we know that $f(0) = 0$.  As a homeomorphism of $V$, $f$ either preserves or reverses the orientation of $V$ (at $0$). (See [1], Section 3.3, for example, for a discussion of orientation on manifolds.) In either case, $f \circ f$ must \textit{preserve} the orientation of $V$ (at $0$).  However, the inverse map $\iota : V \rightarrow V$, which is the linear transformation with matrix $-I_{2n+1}$, reverses orientation since $\det(-I_{2n+1}) = -1$.  Therefore, $f \circ f = \iota$ is not possible.  In summary:

\begin{thm}  An $n$-dimensional vector space $V$ over $\R$ admits a continuous inverse ambiguous function if and only if $n$ is even.
\end{thm}

Since the additive matrix groups $M_n(\R)$ and $M_n(\C)$ are isomorphic as Lie groups to, respectively, $\R^{n^2}$ and $\C^{n^2}$, we immediately obtain the following result.

\begin{cor}  $ $
\begin{enumerate}
    \item For every positive integer $n$, the topological group $M_n(\C)$ admits a continuous inverse ambiguous function, one example being $A \mapsto iA$.
    \item The topological group $M_n(\R)$ admits a continuous inverse ambiguous function if and only if $n$ even.
\end{enumerate} 
\end{cor}

\begin{cor} \label{trace0} $ $
\begin{enumerate}
    \item For every positive integer $n$, the subspace of all matrices in $M_n(\C)$ having trace $0$ admits a continuous inverse ambiguous function.
    \item The subspace of all matrices in $M_n(\R)$ having trace $0$ admits a continuous inverse ambiguous function if and only if $n$ is $odd$.
\end{enumerate}
\end{cor}

\begin{proof}  This follows from the previous corollary in light of the fact that the subspace of trace-0 matrices has codimension 1.
\end{proof}

Consideration of orientation can also be used to determine existence or non-existence of continuous inverse ambiguous functions on $n$-spheres.  The only spheres that are connected Lie groups are $S^1$ and $S^3$ (see [1], Section 3.C).  We saw earlier that there are no continuous ambiguous functions defined on $S^1$ (Theorem~\ref{circleno}).  It turns out the same holds for $S^3$.  The multiplicative group $S^3$ can be identified with the unit quaterions $\R \oplus \R i \oplus \R j \oplus \R k$, where $i^2 = j^2 = k^2 = 1$ and $ijk = -1$.  Now, if $f: S^3 \rightarrow S^3$ were continuous and inverse ambiguous, then $f \circ f$ would preserve the orientation of $S^3$.  However, the group inverse map $\iota$ fixes $1$ and maps each of $i$, $j$, and $k$ to their additive inverses $-i$, $-j$, and $-k$.  Therefore, with respect to the ordered basis $\{1, i, j, k\}$, $\iota$ is given by the matrix

\[\left(\begin{array}{cccc} 1 & 0 & 0 & 0 \\
                           0 & -1 & 0 & 0 \\
                           0 & 0 & -1 & 0 \\
                           0 & 0 & 0 & -1 
       \end{array} \right), \]
\vspace{5pt}

\noindent which has determinant -1.  Therefore, $\iota$ is orientation reversing, and this implies that $f \circ f = \iota$ is not possible.  Hence:

\begin{thm}  No continuous inverse ambiguous functions are defined on $S^3$.
\end{thm}

Finally, even though $S^n$ does not have a group structure for all $n$, the symmetry of each sphere (as a subspace of Euclidean space) allows us to consider homeomorphisms $f : S^n \rightarrow S^n$ such that 
\begin{equation} \label{addinvIA} f(f(z)) = -z \:\:\:\:\:\:\: \textrm{ for all } z \in S^n.
\end{equation}  
(In other words, $f \circ f$ is the antipodal map on $S^n$.) The linear map $\iota : z \mapsto -z$ on $\R^{n+1}$ has matrix $-I_{n+1}$.  If $n$ even, then $\det(-I_{n+1}) = -1$, which indicates that $\iota$ is orientation reversing on $\R^{n+1}$, and thus also on $S^n$.  As a result, $f \circ f$ cannot equal $\iota$ for any homeomorphism of $S^n$ to itself since, as noted before, the composition $f \circ f$ would be orientation preserving regardless of whether $f$ is orientation preserving or reversing.  On the other hand, if $m = 2n - 1$ is odd, it is easily verified that the function defined in Equation (\ref{imap}) restricts to a homeomophism from $S^{m}$ to itself that satisfies Equation (\ref{addinvIA}).  This completes the proof of the following:

\begin{thm}
There exists a homeomorphism $f: S^n \rightarrow S^n$ satisfying the equation $f(f(z)) = -z$ for all $z \in S^n$ if and only if $n$ is odd.
\end{thm}

\section{Multiplicative Matrix Groups}

In this final section, we begin an investigation of the multiplicative group $G = GL_n(\R)$ and some of its subgroups.  We observe that $G$ has two connected components, namely $GL_n^+(\R)$, the normal subgroup of all matrices in $G$ with a positive determinant, and $GL_n^-(\R)$, the subset of all matrices in $G$ with a negative determinant. 
Moreover, the orthogonal group  
\[ O(n) = \{ A \in G \: | \: AA^T = I_n \} \]
is a closed subgroup of $G$ whose components are $SO(n) = O^+(n)$ and its complement.

When $n$ is odd, we can use the disconnectedness of $G$ to define a continuous inverse ambiguous function on $G$.

\begin{thm} \label{pingpong} Assume $n$ is positive odd integer. Then the function $f: GL_n(\R) \rightarrow GL_n(\R)$ defined by
\[ f(A) = \left\{ \begin{array}{lcl}
                -A & & \textrm{if } A \in GL_n^+(\R) \\
                -A^{-1} & & \textrm{if } A \in GL_n^-(\R)
                  \end{array} \right. \]
is continuous and inverse ambiguous.  Furthermore, this function restricts to a continuous inverse ambiguous function from $O(n)$ to $O(n)$.
\end{thm}

\begin{proof}  Suppose $n$ is odd.  First, we see that if $\det(A) > 0$, then \[\det(-A) = (-1)^n\det(A) = -\det(A) < 0.\]  Similarly, if $\det(A) < 0$, then \[\det(-A^{-1}) = (-1)^n\det(A^{-1}) = -\det(A) > 0.\]  This shows that $f$ is well-defined.  The continuity of $f$ is clear, as is the fact that $f(f(A)) = A^{-1}$ for all $A$.  The last statement follows from the fact that $A \in O(n)$ if and only if $-A$ and $A^{-1}$ belong to $O(n)$.
\end{proof}

If $n$ is even, then the strategy used to define the function in the above theorem does not work since $\det(-A) = (-1)^n\det(A) = \det(A)$.  In other words, $A$ and $-A$ belong to the same connected component.  In fact, as we will demonstrate, when $n = 2$, the strategy of defining a continuous inverse ambiguous function that swaps components of $GL_2(\R)$ is not possible.

\begin{lem} In $GL_2^+(\R)$, there are exactly 2 self-invertible elements, while in $GL_2^-(\R)$ there are infinitely many self-invertible elements.
\end{lem}

\begin{proof}  Suppose $A = \left(\begin{array}{cc} a & b \\ c & d \end{array} \right)$ satisfies $A^2 = I_2 = \left(\begin{array}{cc} 1 & 0 \\ 0 & 1 \end{array} \right)$.  Then
\[ a^2 = 1 - bc = d^2, \:\:\:\:\:  \:\:\:\:\: \textrm{and} \:\:\:\:\: \:\:\:\:\: b(a+d) = 0 = c(a+d) \]
are necessary and sufficient conditions for the entries of $A$.

If $a, d \ne 0$, then the first equation gives us $a = \pm d$.  If, further, $a = d$, then $a + d \ne 0$, which implies that $b = c = 0$, and hence $a^2 = 1 - 0 = d^2$, or $a = d = \pm 1$.  The resulting matrices are $\pm I_2$, both of which have determinant 1.  If, $a = -d \ne 0$, then any real numbers $b, c$ satisfying $bc = 1 - a^2$ (for which there are infinitely many combinations for any value of $a$) will give a solution for $A$. In this case, $\det(A) = ad - bc = -a^2 - bc = -1$.  Finally, if $a = d = 0$, then $1 - bc = 0 = ad$, and so $\det(A) = ad - bc = 0 - bc = -1$.
\end{proof}

\begin{rem}  The two self-invertible elements in $GL_2^+(\R)$ are the identity matrix and the matrix representing a half turn about the origin.  Among the self-invertible elements in $GL_2^-(\R)$ are matrices representing reflections across lines through the origin.  Because a self-invertible matrix in $GL_n(\R)$ must have determinant $\pm 1$, the lemma also shows that $SO(2)$ (resp. $O(2) - SO(2)$) has exactly 2 (resp. infinitely many) self-invertible elements.
\end{rem}

In light of this lemma and Theorem~\ref{generalexistence}, we see that any continuous inverse ambiguous function $f : GL_2(\R) \rightarrow GL_2(\R)$ must map $GL_2^+(\R)$ to $GL_2^+(\R)$ and $GL_2^-(\R)$ to $GL_2^-(\R)$.  More specifically, either $f$ fixes $I_2$ and $-I_2$ or transposes them.  Similarly, any continuous inverse ambiguous function $f : O(2) \rightarrow O(2)$ would necessarily restrict to a continuous inverse ambiguous function on $SO(2)$.  However, $SO(2)$ is isomorphic to $S^1$ as Lie groups, and $S^1$ admits no continuous inverse ambiguous functions (Lemma~\ref{circleno}).  This proves the following:

\begin{thm} There are no continuous inverse ambiguous functions defined on $O(2)$ or on $SO(2)$.
\end{thm}

We next consider the special linear group $SL_n(\R)$, which is the Lie subgroup of $GL_n(\R)$ consisting of all matrices with determinant equal to $1$.  Our investigation will consider orientations as we did in the previous section.  We first establish a lemma.

\begin{lem}  If $\iota : GL_n(\R) \rightarrow GL_n(\R)$ is the inverse map, then the differential of $\iota$ at the identity matrix $I_n$ is $-I^{n^2}$.
\end{lem}

\begin{proof} For each $1 \le i, j \le n$, denote by $e_{ij}$ the $n \times n$ matrix with all $0$ entries except for a $1$ in the $ij^\textrm{th}$ position.  Consider the $n^2$ curves in $GL_n(\R)$, defined for $t \in \R$ and each passing through $I_n$ at $t = 0$:
\[ \gamma_{ij}(t) = I_n + t \cdot e_{ij}.\]
One can check that 
$$ \iota(\gamma_{ij}(t)) = \left\{\begin{array}{lcl} 
I_n - t \cdot e_{ij} & & \textrm{if } i \ne j \\ \\
I_n - \dfrac{t}{t+1} \cdot e_{ii}  & & \textrm{if } i = j. \end{array} \right.$$
Differentiating each of these image curves and evaluating at $t = 0$, we get
$d\iota|_{I_n} (e_{ij}) = -e_{ij}$ for all $i, j$.  This implies that $d\iota|_{I_n} = -I_{n^2}$, as desired.
\end{proof}

Now, the determinant of $d\iota|_{I_n}$ is $(-1)^{n^2}$, which indicates that $\iota$ is orientation preserving (resp. reversing) on $GL_n^+(\R)$ when $n$ is even (resp. odd).  Thus, we do \textit{not} arrive at the usual orientation-related contradiction that would rule out the existence of a continuous inverse ambiguous function on $GL^+_n(\R)$ when $n$ is even.  However, we can deduce the following:

\begin{thm}  There are no continuous inverse ambiguous functions on $GL_n^+(\R)$ when $n$ is odd.
\end{thm}

Consequently, the only continuous inverse ambiguous functions on $GL_n(\R)$ for odd $n$ must map each of its two connected components to the other one (as demonstrated in Theorem~\ref{pingpong}).

Now consider the connected, $n^2 - 1$ dimensional, Lie subgroup $SL_n(\R)$ of $GL_n(\R)$.  (See [1], Example 5.27.)  Because $SL_n(\R)$ is an embedded submanifold of the general linear group, $d\iota|_{I_n}$ on the tangent space of $SL_n(\R)$ is simply the restriction of $d\iota|_{I_n}$ from the tangent space of $GL_n(\R)$ [1, Chapter 5].  In other words, for $SL_n(\R)$, $d\iota|_{I^n} = -I_{n^2 - 1}$, which has determinant $-1$ if and only if $n$ is even.  Hence, whether $f$ is a orientation preserving or reversing continuous inverse ambiguous function on $SL_n(\R)$ ($n$ even), then $f \circ f$ is orientation preserving but $\iota$ is orientation reversing.  Therefore,

\begin{thm} If $n$ is even, there are no continuous inverse ambiguous functions on $SL_n(\R)$.
\end{thm}

Likewise, since $SO(n)$ is a connected, $\binom{n}{2}$-dimensional embedded Lie subgroup of $GL_n(\R)$ [1, Example 5.42], then a similar argument shows that whenever $\binom{n}{2} = \frac{1}{2} (n)(n-1)$ is odd, then there are no continuous inverse ambiguous functions defined on $SO(n)$.  More simply,

\begin{thm} If $n \equiv 2$ or $3 \mod 4$, then there are no continuous inverse ambiguous functions on $SO(n)$.
\end{thm}

When the group inverse map is orientation preserving and the matrix group is connected, our methods do not lead to the non-existence of continuous inverse ambiguous functions.  If further results are to be obtained for groups like $GL_n^+(\R)$ ($n$ even), $SL_n(\R)$ ($n$ odd), $SO(n)$ ($n$ congruent to $0$ or $1$ modulo $4$), and $GL_n(\C)$ (any $n > 1$), other tools will need to be developed.

\end{document}